\documentclass{elsart3-1}

\usepackage{amsfonts}
\usepackage{mathrsfs}
\usepackage[dvips]{xcolor}

\usepackage{amssymb}
\usepackage{}

\usepackage{amsmath,amssymb}

\usepackage[english,francais]{babel}
\definecolor{red}{rgb}{1,0,0}

\def\small{\smallskip\noindent}

\def\today{\noindent\number\day
\space\ifcase\month\or
  January\or February\or March\or April\or May\or June\or
  July\or August\or September\or October\or November\or December\fi
  \space\number\year}

\def\dint{\int\kern-.6em\int}

\def\sup {\mathop{{\rm sup}}}

\def\=d{{\,\buildrel (d) \over =\,}}
\def\a.s.{{\buildrel a.s. \over \longrightarrow}}

\newtheorem{theorem}{Theorem}[section]
\newtheorem{lemma}[theorem]{Lemma}
\newtheorem{e-proposition}[theorem]{Proposition}
\newtheorem{corollary}[theorem]{Corollary}
\newtheorem{e-definition}[theorem]{Definition}
\newtheorem{remark}{Remark}
\newtheorem{example}{Example}
\newtheorem{theoreme}{Th\'eor\`eme}[section]

\newtheorem{proposition}[theoreme]{Proposition}

\renewcommand{\theequation}{\arabic{equation}}
\def\phi{\varphi}

\def\proof{{\medskip\noindent {\sl Proof. }}}

\def\ignore#1{}

\def\og{\leavevmode\raise.3ex\hbox{$\scriptscriptstyle\langle\!\langle$~}}
\def\fg{\leavevmode\raise.3ex\hbox{~$\!\scriptscriptstyle\,\rangle\!\rangle$}}
\begin{document}
\centerline{}
\begin{frontmatter}


\selectlanguage{english}
\title{On representation theorem of sublinear expectation related to $G$-L\'{e}vy process and paths of $G$-L\'{e}vy process }


\author{ Liying REN \thanksref{label1}}
\ead{ly\_ren@hotmail.com}
\thanks[label1]{The author thanks the partial support from The National Basic Research Program of
China (973 Program) grant No. 2007CB814900 (Financial Risk).}
\address{School of Mathematics \  Shandong University 250100, Jinan, China}


\begin{abstract}

In this paper, we are concerned with the representation of an
important sublinear expectation $\mathbb{E}^{G}[\cdot]$ under which
framework a new stochastic process $G$-L\'{e}vy process has been
introduced. We show the existence of a weakly compact family of
probability measures $\mathscr{P}$ to give the representation by
using two different methods.

\end{abstract}


\end{frontmatter}

\section{Introduction}

 In 2006, \cite{Peng 2006} introduced the notion of $G$-Brownian motion and $G$-expectation by using a nonlinear partial differential equation called $G$-heat equation.
 Under this fully nonlinear
$G$-expectation which serves as a tool to measure risk and model
uncertainty, corresponding theory of stochastic calculus is further
developed, such as a new type of It\^{o} formula, the existence and
uniqueness for stochastic differential equation driven by
$G$-Brownian motion and martingale representation theorem for
$G$-expectation, etc.(see \cite{Peng 2007a}, \cite{Peng 2008a}, \cite{STZ}).

In \cite{DHP} and \cite{HP 2009b}, it has been proved that a $G$-expectation admits a representation with respect to a
weakly compact set of probability measures.
Recently, Hu and Peng have introduced in \cite{HP 2009a} a new
stochastic process called $G$-L\'{e}vy process
 under a framework of sublinear expectation $\mathbb{E}^{G}$ which derives from an integro-partial differential equation.
In this paper, motivated by both two methods established
respectively in \cite{HP 2009b} and Section 4.1 of \cite{DHP}, we
use both the elementary representation of sublinear expectation and
the stochastic control method to find a weakly relatively compact
family of probability measures 
$\{P_{\theta}: \theta\in \Theta\}$
which gives the following representation:
$$
\mathbb{E}^{G}[X]=\sup
_{\theta\in\Theta}\int_{\Omega}X(\omega)dP_{\theta}(\omega),
$$
where the state space $\Omega$ in our study is $D([0,\infty), \mathbb{R}^{d})$. \\

This paper is organized as follows: in Section 2, we review some
important notions and results of sublinear expectations,
$G$-L\'{e}vy process, capacity  and function space related to an
upper-expectation.  In Section 3, we use the elementary
representation of sublinear expectation by a  family of finitely
additive linear expectations $\{E_{\lambda}:\lambda\in\Lambda\}$ to
find the desired weakly relatively compact family of probability
measures $\{P_{\theta}: \theta\in \Theta\}$. A concrete weakly
relatively compact family of probability measures is constructed
through the method of optimal stochastic control in Section 4. In
Section 5, we give some elementary characterizations of
$L_{ip}(\Omega)$. A generalized Kolmogorov-Chentsov's criterion for
weak relative compactness  has also been presented in the Appendix.

\section{Basic settings }
For a given positive integer $n$, we will denote by
$\langle x,y\rangle$ the scalar product of $x, y\in\mathbb{R}^n$ and
by $|x|=\langle x,x\rangle^\frac{1}{2}$ the Euclidean norm of $x.$
\\

Let $\Omega$ be a given set and let $\mathcal{H}$ be a linear space of
real-valued functions defined on $\Omega$ such that for each integer
$n$, if $X_1, X_2, \cdots, X_n\in\mathcal{H},$ then $\varphi(X_1,
X_2, \cdots, X_n)\in\mathcal{H}$ for each $\varphi\in
C_{b,Lip}(\mathbb{R}^n),$ where $
C_{b,Lip}(\mathbb{R}^n)$ denotes the space of bounded and Lipschitz functions on $\mathbb{R}^n$.
$\mathcal{H}$ is usually considered as a
space of random variables.\\

In the following we will present some preliminaries in the theory of sublinear
expectation and the related $G$-L\'{e}vy process. More details of this section can
be found in \cite{HP 2009a}.

\begin{e-definition}\label{sublinear}
A functional $\mathbb{\hat E}: \mathcal{H} \rightarrow \mathbb{R}$
is called a sublinear expectation if it satisfies the following
properties: for all $X, Y \in\mathcal{H},$ we have\\
(a) Monotonicity: if $X\geq Y,$ then $\mathbb{\hat E}[X]\geq
\mathbb{\hat E}[Y].$\\
(b) Constant preserving: $\mathbb{\hat E}[c]= c, \forall
c\in\mathbb{R}.$\\
(c) Sub-additivity: $\mathbb{\hat E}[X]-\mathbb{\hat E}[Y]\leq
\mathbb{\hat E}[X-Y].$\\
(d) Positive homogeneity: $\mathbb{\hat E}[\lambda
X]=\lambda\mathbb{\hat E}[X], \forall \lambda \geq0.$
\\
We call the triple $(\Omega, \mathcal{H}, \mathbb{\hat E})$ a sublinear
expectation space similar to the probability space $(\Omega,
\mathcal{F}, P).$
\end{e-definition}

\begin{e-definition}
 Let $X_{1}$ and $X_{2}$ be two n-dimensional random vectors defined respectively on the sublinear expectation spaces
$(\Omega_{1}, \mathcal{H}_{1}, \mathbb{\hat E}_{1})$ and $(\Omega_{2}, \mathcal{H}_{2}, \mathbb{\hat E}_{2}).$ They are called identically
distributed if  $\mathbb{\hat E}_{1}[\varphi(X_{1})]=\mathbb{\hat
E}_{2}[\varphi(X_{2})]$ for all $\varphi\in C_{b, Lip}(\mathbb{R}^n)$, denoted by $X_{1}\stackrel{D}{=}X_{2}.$
\end{e-definition}

\begin{e-definition}
On a sublinear expectation space $(\Omega, \mathcal{H}, \mathbb{\hat
E}),$ an n-dimensional random vector $Y$ is said to be independent
of another m-dimensional random vector $X$, if for each test
function $\varphi\in C_{b,Lip}(\mathbb{R}^{m}\times\mathbb{R}^{n}),$ we have
$$
\mathbb{\hat E}[\varphi(X, Y)]=\mathbb{\hat E}[\mathbb{\hat
E}[\varphi(x, Y)]_{x=X}].
$$
\end{e-definition}

\subsection {$G$-L\'{e}vy process and Sublinear expectation related to $G$-L\'{e}vy process}
\begin{e-definition} \label{levy proc}{\bf (L\'{e}vy process)}
Let $X=(X_{t})_{t\geq0}$ be a d-dimensional c\`{a}dl\`{a}g process  defined on a sublinear expectation space $(\Omega, \mathcal{H}, \mathbb{\hat
E}).$ We say that $X$ is a L\'{e}vy process if:
\\
(a) $X_{0}=0;$ \\
(b) Independent increments: for each $t,s>0,$ the increments $X_{t+s}-X_{t}$ is independent
of $(X_{t_{1}}, X_{t_{2}}, \ldots , X_{t_{n}}),$ for each $n \in \mathbb{N}$ and $0\leq t_{1}\leq\ldots\leq t_{n}\leq t;$\\
(c) Stationary increments: the distribution of $X_{t+s}-X_{t}$ does not depend on t.\\
\end{e-definition}

\begin{e-definition}{\bf ($G$-L\'{e}vy process)}
A d-dimensional process $(X_{t})_{t\geq0}$ on a sublinear expectation space $(\Omega, \mathcal{H}, \mathbb{\hat
E})$ is called a G-L\'{e}vy process if:\\
(a) $X$ satisfies (a)-(c) in Definition. \ref{levy proc};\\
(b) for each $t\geq0$, there exists a decomposition $X_{t}=X_{t}^{c}+X_{t}^{d};$ \\
(c) $(X_{t}^{c},X_{t}^{d})_{t\geq0}$ is a 2d-dimensional L\'{e}vy process satisfying
\begin{equation}\label{assup2}
\lim _{t\downarrow0}\mathbb{\hat E}[|X_{t}^{c}|^{3}]t^{-1}=0; \
\mathbb{\hat E}[|X_{t}^{d}|]\leq Ct \ {\rm for} \  t\geq0,
\end{equation}
where C is a constant.
\end{e-definition}

\begin{remark}
We know that in classical linear expectation case, assumption (b) of $(X_{t})_{t\geq0}$ obviously holds by the L\'{e}vy-It\^{o} decomposition, where $(X_{t}^{c})_{t\geq0}$
and $(X_{t}^{d})_{t\geq0}$ represent the continuous part and jump part respectively.
By assumption (\ref{assup2}) on $(X_{t}^{c})_{t\geq0}$ and $(X_{t}^{d})_{t\geq0}$,
we know that $(X_{t}^{c})_{t\geq0}$ is a generalized G-Brownian motion, and $(X_{t}^{d})_{t\geq0}$ is of finite variation.
\end{remark}

\begin{theorem} Let $(X_{t})_{t\geq0}$ be a $d$-dimensional $G$-L\'{e}vy process with the decomposition $X_{t}=X_{t}^{c}+X_{t}^{d}.$ For each $\varphi \in C_{b,Lip}(\mathbb{R}^n),$
$u(t,x):=\mathbb{\hat E}[\varphi(x+X_{t})]$ is a viscosity solution of the following
equation:
$$
\partial_{t}u(t,x)-G_{X}[u(t,x+\cdot)-u(t,x)]=0, u(0,x)=\varphi(x),
$$
 where $G_{X}[f(\cdot)]$ is a nonlocal operator defined by
 $$
 G_{X}[f(\cdot)]:=\lim _{\delta\downarrow 0}\mathbb{\hat E}[f(X_{\delta})]\delta^{-1}
 \ {\rm for} \ \  f \in C_{b}^{3}(\mathbb{R}^{d}) \ {\rm with} \  f(0)=0.
 $$
\end{theorem}

Now we will consider a particular sublinear
expectation related to $G$-L\'{e}vy process, and throughout this paper our study
will be worked under such a sublinear expectation.
Let $\Omega=D_0([0,\infty),\mathbb{R}^{d})$ denote the space of all
$\mathbb{R}^{d}-$valued c\`{a}dl\`{a}g paths $(\omega_t)_{t\geq0}$
with $\omega_0=0,$ equipped with the distance (introduced in \cite{Bill})
$$
d_{\infty}^{\circ}(\omega^1,
\omega^2):=\sum_{i=1}^{\infty}2^{-i}[d_{i}^{\circ}(\omega^1_i, \omega_i^2)\wedge1],
$$
 where $d_{t}^{\circ}$ is a metric defined on $D_0([0,t],\mathbb{R}^{d})$. $D_0([0,t],\mathbb{R}^{d})$ is complete and separable under $d_{t}^{\circ}$, hence is a Polish space. $D_0([0,\infty),\mathbb{R}^{d})$ is also a Polish space under the Skorohod topology $d_{\infty}^{\circ}$. Let $\mathscr{B}(\Omega)$ denote the $\sigma$-algebra generated by
all open sets. We will consider the canonical process $B_t(\omega)=\omega_{t}$ for
$\omega \in \Omega,\  t\geq0.$ We introduce the space of finite dimensional cylinder random variables: for each fixed $T\geq0,$ we denote $\Omega_{T}=\{\omega_{\cdot\wedge T } : \omega\in \Omega\}$ and set
\par


%

$$\begin{array}{ccl}
L_{ip}(\Omega_T)
 &:= &\{\varphi(B_{t_1\wedge T}, B_{t_2\wedge T},
\cdots, B_{t_n\wedge T}):  n\in\mathbb{N}, t_1, t_2, \cdots,
t_n\in[0, \infty), \varphi\in C_{b,Lip}(\mathbb{R}^{d\times n})\}; \\[10pt]
 L_{ip}(\Omega)& := & \{\varphi(B_{t_1}, B_{t_2}, \cdots,
B_{t_n}):  n\in\mathbb{N}, t_1, t_2, \cdots, t_n\in[0, \infty), \varphi\in
C_{b,Lip}(\mathbb{R}^{d\times n})\}.
\end{array}
$$
It is clear that for $t\leq T,$ $ L_{ip}(\Omega_t)\subseteq
L_{ip}(\Omega_T)$  and $L_{ip}(\Omega)=
\cup_{n=1}^{\infty}L_{ip}(\Omega_n).$


%

A sublinear expectation $\mathbb{\hat E}[\cdot]$ defined on $L_{ip}(\Omega)$ through the following procedure is a sublinear expectation related to the $G$-L\'{e}vy process:\\

Step 1. For each $\xi\in L_{ip}(\Omega)$
with  the form $\xi=\varphi(B_{t+s}-B_{t}),$ $ t,s\geq 0$ and $\varphi \in C_{b,Lip}(\mathbb{R}^n).$
We define $\mathbb{\hat E}[\xi]=u(s,0),$
where u is
a viscosity solution of the following integro-partial differential equation:
\begin{equation}\label{ipde}
\dfrac{\partial u}{\partial t}
 = G(u,Du, D^{2}u), \ \  on \ (t,x)\in [0,T)\times \mathbb{R}^{d}
 \end{equation}
with the initial condition $u(0,x)=\varphi(x),$
where
$$G(u,p,A)=\sup _{(\nu,q,Q)\in \mathcal{U}}\{\int_{\mathbb{R}^{d}\backslash \{0\}}(u(t,x+z)-u(t,x))\nu(dz)+\langle p,q\rangle+ \frac{1}{2}tr[AQQ^{T}]\},$$ and $\mathcal{U}$ represents $G_{X}$.\\

Step 2. For each $\xi \in L_{ip}(\Omega),$ there exists a $\psi \in C_{b,Lip}(\mathbb{R}^{d\times m})$ such that $\xi=\psi(B_{t_1}-B_{t_0}, B_{t_2}-B_{t_1}, \cdots,
B_{t_m}-B_{t_{m-1}})$  for some $t_{1}<t_{2}<\ldots<t_{m}.$ Then $\mathbb{\hat E}[\xi]$ is defined by $\psi_{m}$ via the following procedure:

$$\begin{array}{rcl}
\psi_{1}(x_{1},\ldots,x_{m-1}) & = & \hat{\mathbb{E}}(\psi(x_{1},\ldots,x_{m-1},B_{t_{m}}-B_{t_{m-1}}));\\[10pt]
\psi_{2}(x_{1},\ldots,x_{m-2}) & = & \hat{\mathbb{E}}(\psi_{1}(x_{1},\ldots,x_{m-2},B_{t_{m-1}}-B_{t_{m-2}}));\\[10pt]
\vdots\\
\psi_{m-1}(x_{1}) & = & \hat{\mathbb{E}}(\psi_{m-2}(x_{1},B_{t_{2}}-B_{t_{1}}));\\[10pt]
\psi_{m} & = & \hat{\mathbb{E}}(\psi_{m-1}(B_{t_{1}})).\end{array}$$
The related conditional expectation of \ $\xi$ under $\Omega_{t_{i}}$ is defined by
$$\begin{array}{ccl}
\hat{\mathbb{E}}[\xi|\Omega_{t_{j}}] & = & \hat{\mathbb{E}}[\psi(B_{t_{1}},B_{t_{2}}-B_{t_{1}},\ldots,B_{t_{m}}-B_{t_{m-1}})|\Omega_{t_{j}}]\\[10pt]
 & = & \psi_{m-j}(B_{t_{1}},\ldots,B_{t_{j}}-B_{t_{j-1}}).\end{array}$$
It can be proved that the canonical process $(B_{t})_{t\geq 0}$ is a $G$-L\'{e}vy process, and $\hat{\mathbb{E}}[\cdot]$ consistently defines a sublinear expectation on $L_{ip}(\Omega)$.
Then for $p\geq1$, the topological completion of $L_{ip}(\Omega_{T})$(resp. $L_{ip}(\Omega)$ under the Banach norm $\parallel\cdot\parallel_p:=\mathbb{\hat E}[|\cdot|^p]^\frac{1}{p}$  is denoted by  $L_G^p(\Omega_{T})$ (resp.
$L_G^p(\Omega)$).
$\hat{\mathbb{E}}[\cdot]$ can be extended continuously and uniquely from $L_{ip}(\Omega_{T})$
(resp. $L_{ip}(\Omega)$) into $L_G^p(\Omega_{T})$ (resp.
$L_G^p(\Omega)$), and it is still a sublinear expectation on
the corresponding completed space. \\

\subsection{Capacity Associated to an Upper Probability and related functional spaces}
Let $\Omega$ be a complete separable metric space equipped with the distance $d$,
$\mathscr{B}(\Omega)$ the Borel $\sigma$-algebra of $\Omega$ and $\mathcal {M}$
the collection of all probability measures on $(\Omega,\mathscr{B}(\Omega)).$
We denote by $L^{0}(\Omega)$ the space of all $\mathscr{B}(\Omega)$- measurable real functions and $C_{b}(\Omega)$ all bounded continuous functions.
For a given subset $\mathscr{P}\subseteq \mathcal {M},$ we denote
$$
c(A):=\sup _{P\in \mathscr{P}}P(A), \ \ A\in \mathscr{B}(\Omega).
$$
It is easy to verify that $c(\cdot)$  satisfies the following theorem.
\begin{theorem}
The set function $c(\cdot)$ is a Choquet capacity,  i.e. (see \cite{Choquet}, \cite{Dellacherie}) \\
(a) $0\leq c(A)\leq 1, \ \  \forall A\subset\Omega.$\\
(b) If $A\subset B,$  then $c(A)\leq c(B).$\\
(c) If $(A_{n})_{n=1}^{\infty}$ is a sequence in $\mathscr{B}(\Omega),$ then $c(\bigcup A_{n})\leq \sum c(A_{n}).$\\
(d) If $(A_{n})_{n=1}^{\infty}$ is an increasing sequence in $\mathscr{B}(\Omega)$: $A_{n}\uparrow A=\bigcup A_{n},$ then $c(\bigcup A_{n})=\lim _{n\rightarrow \infty}c(A_{n}).$\\
\end{theorem}

\begin{e-definition}
We say that a set  A is polar if $c(A)=0$ and a property holds
"quasi-surely" (q.s.) if it holds outside a polar set. In other
words, $A\in \mathscr{B}(\Omega)$ is polar if and  only if $P(A)=0$
for any $P\in \mathscr{P}.$
\end{e-definition}

By applying the Borel-Cantelli Lemma we could get immediately the following lemma.
\begin{lemma}\label{bc}
Let $(A_{n})_{n \in \mathbb{N}}$ be a sequence of Borel sets such that
$$
\sum_{n=1}^{\infty}c(A_{n})<\infty.
$$
Then $\limsup_{n\rightarrow \infty} A_{n}$ is polar.
\end{lemma}
The upper expectation (see \cite{Huber}) of $\mathscr{P}$ is defined as follows:
for each $X\in L^{0}(\Omega)$ such that $E_{P}[X]$ exists for each $P\in \mathscr{P}$,
$$\mathbb{E}[X]=\mathbb{E}^{\mathscr{P}}[X]:=\sup _{P\in \mathscr{P}}E_{P}[X].$$

%

If for $p>0$, we denote
\begin{itemize}
\item $\mathcal {L}^{p}:=\{X\in L^{0}(\Omega): \mathbb{E}[|X|^{p}]=\sup _{p\in \mathscr{P}}E_{P}[|X|^{p}]<\infty\} ;$
\item $\mathcal {N}^{p}:=\{X\in L^{0}(\Omega): \mathbb{E}[|X|^{p}]=0\};$
\item $\mathcal {N}:=\{X\in L^{0}(\Omega): X=0, c-q.s.\}.$
\end{itemize}

It is seen that $\mathcal {L}^{p}$ and $\mathcal {N}^{p}$ are linear spaces and $\mathcal {N}^{p}=\mathcal {N}$, for $p>0$.
All of the following definitions and propositions can be found in \cite{DHP}, and the proofs are omitted here.
\begin{proposition}
For each  $p\geq1$,  $\mathbb{L}^{p}:=\mathcal {L}^{p}/\mathcal {N}$ is a Banach space under the
norm $\|X\|_{p}:=\big( \mathbb{E} [|X|^{p}] \big)^{1/p}$; for each
$p<1$, $\mathbb{L}^{p}$ is a complete metric space under the distance
$d(X,Y):= \mathbb{E} [|X-Y|^{p}]$.
\end{proposition}

With respect to the distance defined on $\mathbb{L}^{p}$, $p>0$, we denote:

\begin{itemize}
\item $\mathbb{L}_{c}^{p}$ the completion of $C_{b}(\Omega)$;
\item $L_{G}^{p}(\Omega)$ the completion of $L_{ip}(\Omega)$.
\end{itemize}

%
\begin{e-definition}
A mapping $X$ on $\Omega$ with values in a topological space is said to be quasi-continuous if
for any $\varepsilon>0$, there exists an open subset $O$ with $c(O)<\varepsilon$ such that $X|_{O^{c}}$ is continuous.
\end{e-definition}


\begin{e-definition}
We say that $X:\Omega\rightarrow\mathbb{R}$ has a quasi-continuous version if there exists a quasi-continuous function $Y:\Omega\rightarrow\mathbb{R}$ with $X=Y$ q.s..
\end{e-definition}
\begin{proposition}
For each  $p>0$,
$$\mathbb{L}^{p}_{c}=\{X\in \mathbb{L}^{p}: X \ has \  a \ quasi-continuous \ version, \lim _{n\rightarrow \infty} \mathbb{E}[|X|^{p}I_{\{|X|>n\}}]=0 \}.$$
\end{proposition}
\begin{proposition}\label{wc}
Let $\{P_{n}\}_{n=1}^{\infty}\subset\mathscr{P}$ converge weakly to $P\in \mathscr{P}$. Then
for each $X\in\mathbb{L}_{c}^{1}$, we have $E_{P_{n}}[X]\rightarrow E_{P}[X]$ as $n\rightarrow\infty$.
\end{proposition}
\section{Representation of the sublinear expectation related to $G$-L\'{e}vy process as an upper-Expectation}
Let $\bar{\Omega}=(\mathbb{R}^{d})^{[0,\infty)}$ denote the space of all
$\mathbb{R}^{d}$-valued functions $(\bar{\omega}_{t})_{t\geq0}$
and $\mathscr{B}(\bar{\Omega})$ denote the $\sigma$-algebra generated by all
finite dimensional cylinder sets. The corresponding canonical process
$\bar{B}_t(\bar{\omega})=\bar{\omega}_{t}$ for
$\bar{\omega} \in \bar{\Omega},\  t\geq0.$  The space of Lipschitizian cylinder
functions on $\bar{\Omega}_{T}$ is denoted by

$$
L_{ip}(\bar{\Omega}):=\{\varphi(\bar{B}_{t_1}, \bar{B}_{t_2}, \cdots,
\bar{B}_{t_n}):  n\in\mathbb{N}, t_1, t_2, \cdots, t_n\in[0, \infty), \varphi\in
C_{b,Lip}(\mathbb{R}^{d\times n})\}
$$
and
$$
L_{ip}(\bar{\Omega}_T):=\{\varphi(\bar{B}_{t_1\wedge T}, \bar{B}_{t_2\wedge T}, \cdots,
\bar{B}_{t_n\wedge T}):  n\in\mathbb{N}, t_1, t_2, \cdots, t_n\in[0, \infty), \varphi\in
C_{b,Lip}(\mathbb{R}^{d\times n})\}.
$$
%
Following the same procedure as the construction of $\hat{\mathbb{E}}[\cdot]$, we can also construct a sublinear expectation $\bar{\mathbb{E}}$ on $(\bar{\Omega},L_{ip}(\bar{\Omega}))$ such that
$(\bar{B}_{t}(\bar{\omega}))_{t\geq0}$ is also a $G$-L\'{e}vy process.\\

The following lemmas can be found in \cite{HP 2009b} and \cite{Peng 2008b}.

\begin{lemma}\label{lem1}
Let $\hat{\mathbb{E}}$ be a sublinear  functional defined on a linear space $\mathcal {H}$,
i.e. (c) and (d) of Definition \ref{sublinear} hold for $\hat{\mathbb{E}}$. Then there exists
a family $\mathscr{Q}=\{E_{\theta}: \theta\in\Theta\}$ of linear functionals defined on $\mathcal {H}$ such that
$$
\hat{\mathbb{E}}[X]:= \sup_{\theta\in\Theta} E_{\theta}[X], \ for \ X\in\mathcal {H}.
$$
and such that, for each $X\in\mathcal{H}$, there exists a $\theta\in\Theta$ such that $\hat{\mathbb{E}}[X]= E_{\theta}[X].$ Moreover if $\hat{\mathbb{E}}$ is a sublinear functional
defined on a linear space $\mathcal {H}$ of functions on $\Omega$ such that (a) of Definition
\ref{sublinear} holds (resp. (a),(b) hold) for $\hat{\mathbb{E}},$
 then (a) also holds (resp. (a),(b) hold) for
 $E_{\theta}, \  \theta\in\Theta.$
\end{lemma}

\begin{lemma}
Let $0\leq t_{1} <t_{2}<\cdots <t_{m}<\infty$ and $\{\varphi_{n}\}_{n=1}^{\infty}\subset \mathcal{C}_{b,lip}(\mathbb{R}^{d\times m})$ satisfy $\varphi_{n}\downarrow 0,$ then $\bar{\mathbb{E}}[\varphi_{n}(\bar{B}_{t_1}, \bar{B}_{t_2}, \cdots,
\bar{B}_{t_m})]\downarrow 0.$
\end{lemma}

\begin{lemma}\label{lem2}
Let $E$ be a finitely additive linear expectation dominated by $\bar{\mathbb{E}}$ on $L_{ip}(\bar{\Omega}),$ then there exists a unique probability measure $Q$ on $(\bar{\Omega},\mathscr{B}(\bar{\Omega}))$ such that $E[X]=E_{Q}[X]$ for each $X\in L_{ip}(\bar{\Omega}).$
\end{lemma}

\begin{remark}
This is a direct result of  Daniell-Stone's theorem and Kolmogorov's consistent theorem
with the help of above lemma.
\end{remark}
By  Lemma \ref{lem1} and Lemma \ref{lem2}, it is easy to get the next result crucial for our following discussions.
\begin{lemma}\label{lem3}
There exists a family of probability measures $\mathscr{P}_{e}$ on $(\bar{\Omega},\mathscr{B}(\bar{\Omega}))$ such that
$$\bar{\mathbb{E}}[X]=\max_{Q\in \mathscr{P}_{e}} E_{Q}[X],  \forall X\in L_{ip}(\bar{\Omega}).$$
\end{lemma}
We denote the associated capacity to $\mathscr{P}_{e}$ by
$$
\tilde{c}(A):=\sup _{Q\in \mathscr{P}_{e} }Q(A), \ \ A\in \mathscr{B}(\bar{\Omega}).
$$
and related upper expectation for each $\mathscr{B}(\bar{\Omega})$-measurable real function
X which makes the following definition meaningful,
$$
\tilde{\mathbb{E}}[X]:=\sup _{Q\in \mathscr{P}_{e} }E_{Q}[X].
$$

\begin{e-definition}
Let I be a set of indices, $(X_{t})_{t\in I}$ and $(Y_{t})_{t\in I}$ be two processes
indexed by I. Y is said to be a quasi-modification of X in for all $t\in I,$ $X_{t}=Y_{t}$ q.s..
\end{e-definition}

\begin{e-definition}
Let $\varepsilon$ be a positive number. A function $y=X(t)$ is said to have
no fewer than m $\varepsilon$-oscillations on a closed interval $[a,b]$ if there exist points
$a\leq t_{0}<t_{1}<\cdots<t_{m}\leq b,$ such that $|X_{t_{k-1}}-X_{t_{k}}|>\varepsilon$ for
$k=1,2,\cdots, m$.
\end{e-definition}
\begin{remark}
We know that for a function $y=X(t)$ to be c\`{a}dl\`{a}g on a interval $[a,b]$, it is necessary and
sufficient that for arbitrary $\varepsilon>0$, it has only finitely many $\varepsilon$-oscillations on $[a,b]$.
\end{remark}

\begin{theorem}\label{kcc}{\bf (Kolmogorov-Chentsov's Criterion)}
Let $(X_{t})_{t\in [0,1]}$ be a stochastic process such that for all $t\in[0,1]$,
 $X_{t}$ belongs to $\mathbb{L}^{1}$ . If it satisfies the following conditions:\\
(i)for any $ t\in[0,1]$, there exists $a>0$,  such that $\lim_{s\rightarrow t}\mathbb{E}[|X_{s}-X_{t}|^{a}]=0$;\\
(ii)for some $C, r >0, p, q\geq0$ with $p+q>0$,  and all $0\leq s\leq u\leq t \leq 1$, it holds that
\begin{equation}\label{kchenstov}
\mathbb{E}[|X_{t}-X_{u}|^{p}|X_{u}-X_{s}|^{q}]\leq C|t-s|^{1+r},
\end{equation}
then  it admits a c\`{a}dl\`{a}g modification.
\end{theorem}

\proof It follows from condition (i) that, for any $\varepsilon>0$ and $t\in[0,T]$,  $\lim_{s\rightarrow t}c(|X_{s}-X_{t}|\geq\varepsilon)=0$. Now let us take the set $J$ of all dyadic numbers belonging to [0,1],\\
\begin{equation*}
J=\Big\{\frac{i}{2^{n}}; n\in \mathbb{N} , i\in \{ 0,1,\cdots, 2^{n} \} \Big\}.
\end{equation*}
It follows from (\ref{kchenstov}) and Chebyshev's inequality that

\begin{equation}
\begin{array}{ccl}\label{i}
c\big(\{|X_{t}-X_{u}|\geq \epsilon_{1}\} \bigcap \{|X_{u}-X_{s}|\geq \epsilon_{2}\}\big)& \leq & \mathbb{E}\big[\dfrac{|X_{t}-X_{u}|^{p}}{\epsilon_{1}^{p}}
\dfrac{|X_{u}-X_{s}|^{q}}{\epsilon_{2}^{q}}\big]\\[10pt]
 & \leq & \dfrac{C|t-s|^{1+r}}{\epsilon_{1}^{p}\epsilon_{2}^{q}}.\end{array}
\end{equation}
Let $A_{k,n}$ denote the event
$$
\Big\{\big|X\big(\frac{k}{2^{n}}\big)-X\big(\frac{k-1}{2^{n}}\big) \big|<\varepsilon_{n}\Big\},
$$
where
$$\begin{array}{ccl}
\varepsilon_{n} & = & \dfrac{1}{2}C^{1/p+q}\dfrac{r(n-1)}{2(p+q)}=L\alpha^{n},\\[10pt]
L & = & (2^{r/2}C)^{1/p+q}, \ \alpha=2^{-r/2(p+q)}<1,
\end{array}\
$$
and
$$B_{kn}=A_{k,n}\bigcup A_{k+1,n}=\Big\{\big|X\big(\frac{k}{2^{n}}\big)-X\big(\frac{k-1}{2^{n}}\big) \big|<\varepsilon_{n}\Big\}\bigcup\Big\{\big|X\big(\frac{k+1}{2^{n}}\big)-X\big(\frac{k}{2^{n}}\big) \big|<\varepsilon_{n}\Big\}.$$
From the inequality (\ref{i}), we have $$c(B_{kn}^{c})\leq 2^{-(1+r/2)(n-1)}, \  k=1,2,\cdots 2^{n}-1,$$
where $B_{kn}^{c}$ is the complementary event of $B_{kn}$.
Let us define
$$D_{n}:=\bigcap_{m=n}^{\infty}\bigcap_{k=1}^{2^{m}-1}B_{km} \ , \  \  D:=\bigcup_{n=1}^{\infty} D_{n},$$
with the complementary events
$$D_{n}^{c}=\bigcup_{m=n}^{\infty}\bigcup_{k=1}^{2^{m}-1}B_{km}^{c} \ , \  \  D^{c}=\bigcap_{n=1}^{\infty} D_{n}^{c}.$$
Then we  have
$$
c(D_{n}^{c})\leq \sum_{m=n}^{\infty}\sum_{k=1}^{2^{m-1}}c(B_{km}^{c})\leq
\sum_{m=n}^{\infty}\sum_{k=1}^{2^{m-1}}2^{-(1+r/2)(m-1)}=\dfrac{2\beta^{n-1}}{1-\beta},
$$
where $\beta=2^{-r/2}<1$, by Lemma \ref{bc} it follows that $c(D^{c})=0,
$
i.e. $D^{c}$ is a polar set. Let us choose a sample function $X(t)$
for which the event $A_{k+1,n}\bigcap D_{n+1}$ occurs. Thanks to the
classical Kolmogorov-Chentsov's criterion, we can continue to get
the result which is essentially in the spirit of Kolmogorov's
criterion, though its proof is much more difficult. We will only outline the
proof for the convenience of readers. An interested reader may refer
to \cite[pp. 159-164]{GS} for details. Suppose that the event $D$ occurs. Then
beginning with some $n_{0}$, all the events $D_{n}$ for $n\geq
n_{0}$ occur for the sample function of the process. For arbitrary
$\varepsilon>0$ we can find an $n\geq n_{0}$ such that
$2L\alpha^{n}/(1-\alpha)^{2}<\varepsilon.$ On an arbitrary set of
the form
$$
J\bigcap \Big[\frac{k-1}{2^{n}}, \frac{k+1}{2^{n}}\Big],
$$
there can be no more than a single $\varepsilon$-oscillation, which yields that the
function $X(t)$ has no more than $2^{n}$ $\varepsilon$-oscillations on $J$. Thus outside a polar set, the sample functions of process $X(t)$ have only finitely many $\varepsilon$-oscillations; that is, the process $X$ has a c\`{a}dl\`{a}g modification. $\Box$
\begin{remark}
This theorem is a generalized Kolmogorov-Chentsov's criterion for c\`{a}dl\`{a}g modification with respect to capacity.
\end{remark}
\begin{lemma}\label{lem4}
For $\bar{B}=\{\bar{B}_{t}:  t\geq 0\},$ there exists  a c\`{a}dl\`{a}g modification
$\tilde{B}=\{\tilde{B}_{t}:  t\geq 0\}$ of $\bar{B}$ $(i.e. \ \tilde{c}(\{\bar{B}_{t}\neq\tilde{B}_{t} \})=0,$ for each $t\geq0$) such that $\tilde{B}_{0}=0.$
\end{lemma}
\proof  The canonical process $(\bar{B}_{t})$ has a decomposition $\bar{B}_{t}=\bar{B}_{t}^{c}+\bar{B}_{t}^{d}$. The part $\bar{B}_{t}^{c}$ has a continuous modification $\tilde{B}_{t}^{c}$. For the other part $\bar{B}_{t}^{d}$, on the basis of (\ref{assup2}) in
the definition of $G$-L\'{e}vy process and Lemma \ref{lem3}, we have that,
for all $0\leq s\leq u\leq t <\infty,$

$$\tilde{\mathbb{E}}\big[|\bar{B}_{t}^{d}-\bar{B}_{u}^{d}||\bar{B}_{u}^{d}-\bar{B}_{s}^{d}|\big]=
\bar{\mathbb{E}}\big[|\bar{B}_{t}^{d}-\bar{B}_{u}^{d}||\bar{B}_{u}^{d}-\bar{B}_{s}^{d}|\big]\leq  C^{2}|t-s|^{2},$$
where the constant $C$ is a constant. By Theorem \ref{kcc}, there exists a c\`{a}dl\`{a}g modification $\tilde{B}^{d}$ of $\bar{B}^{d}$. Consequently, $\tilde{B}=\{\tilde{B}_{t}: \tilde{B}_{t}=\tilde{B}_{t}^{c}+\tilde{B}_{t}^{d}, t\geq0\}$ is the desired c\`{a}dl\`{a}g modification of $\bar{B}$. $\Box$
\\

\par
The family of probability measures $\mathscr{P}_{e}$ on $(\bar{\Omega}, \mathscr{B}(\bar{\Omega})$ introduces a new family
of probability measures $\mathscr{P}_{1}:=\{Q\circ\tilde{B}^{-1}: Q\in\mathscr{P}_{e}\}$ on $(\Omega, \mathscr{B}(\Omega))$.

\begin{lemma} The family of probability measures $\mathscr{P}_{1}$ is tight.
\end{lemma}
\proof
Since for all $0\leq s\leq u\leq t <\infty,$ there exists $C'>0$, such that
$$\begin{array}{rcl}
\tilde{\mathbb{E}}\big[|\tilde{B}_{t}^{c}-\tilde{B}_{s}^{c}|^{4}\big] & = & \tilde{\mathbb{E}}\big[|\bar{B}_{t}^{c}-\bar{B}_{s}^{c}|^{4}\big]\leq C'|t-s|^{2},\\[10pt]
\tilde{\mathbb{E}}\big[|\tilde{B}_{t}^{d}-\tilde{B}_{u}^{d}||\tilde{B}_{u}^{d}-\tilde{B}_{s}^{d}|\big] & = & \tilde{\mathbb{E}}\big[|\bar{B}_{t}^{d}-\bar{B}_{u}^{d}||\bar{B}_{u}^{d}-\bar{B}_{s}^{d}|\big]\leq C'|t-s|^{2}.
\end{array}\
$$
Due to the generalized Kolmogorov -Chentsov's criterion for tightness (see Theorem \ref{klc}) combined with Corollary \ref{kca}, this implies the wished
statement. $\Box$

By Lemma \ref{lem3} and \ref{lem4}, the representation of sublinear expectation related to $G$-L\'{e}vy process with respect to $\mathscr{P}_{1}$ is given by the following theorem.

\begin{theorem}
For each monotonic and sublinear function $G_{X}[f(\cdot)]: \  \mathbb{R}^{d}\mapsto \mathbb{R},$ where $f\in \  C_{b}^{3}(\mathbb{R}^{d}) $  with $f(0)=0$, let $\mathbb{E}^{G}$ be the corresponding sublinear expectation on $(\Omega, L_{ip}(\Omega))$.
Then there exists a relatively compact family of probability measures $\mathscr{P}_{1}$ on $(\Omega, L_{ip}(\Omega))$, such that
$$
\mathbb{E}^{G}[X]=\max_{P\in\mathscr{P}_{1}} E_{P}[X],  \   \forall X \in L_{ip}(\Omega).
$$
\end{theorem}

\section{Representation of $\mathbb{E}^{G}$ using the stochastic control method}
In this section we will construct a family of probability measures on $\Omega$ for which the upper expectation coincides with the sublinear expectation $\mathbb{E}^{G}$ on $L_{ip}(\Omega)$ through a method of optimal stochastic controls introduced in \cite{DHP}.

Let $(\Omega, \mathcal {F}, P)$ be a probability space and $(L_{t})_{t\geq0}=(L_{t}^{i})_{i=1}^{d}$ a $d$-dimensional L\'{e}vy process with finite variation in this space. The filtration generated by $L$ is denoted by
$$\mathcal {F}_{t}:=\sigma\{L_{u}, 0\leq u\leq t\}\vee\mathcal {N}, \ \ \mathbb{F}:=\{\mathcal {F}_{t}\}_{t\geq0},$$
where $\mathcal {N}$ is the collection of $P$-null subsets. We also denote, for a fixed $s\geq0,$
$$\mathcal {F}_{t}^{s}:=\sigma\{L_{s+u}-L_{s}, 0\leq u\leq t\}\vee\mathcal {N}, \ \ \mathbb{F}^{s}:=\{\mathcal {F}_{t}^{s}\}_{t\geq0}.$$
Let $\Theta=(\Theta^{c},\Theta^{d})=({}^{1}\!\Theta^{c}, {}^{2}\!\Theta^{c}, \Theta^{d})$ be a given  bounded and closed subset in $\mathbb{R}^{d\times 3d}$. We denote a collection of $\Theta$-valued processes on an interval $[t,T]\subset[0,\infty)$ by
$\mathcal {A}_{t,T}^{\Theta}:=\{\theta=(\theta^{c}, \theta^{d})=({}^{1}\!\theta^{c}, {}^{2}\!\theta^{c}, \theta^{d}): \theta^{c}$ is $\Theta^{c}$-valued $\mathbb{F}$-adapted  and $\theta^{d}$ is $\Theta^{d}$-valued $\mathbb{F}$-predictable\}.
For each fixed $\theta=(\theta^{c},\theta^{d}) \in \mathcal {A}_{t,T}^{\Theta}$ we denote
$$B_{T}^{t,\theta}:=\int_{t}^{T}\theta_{s}d(L_{s}^{c},L_{s}^{d})^{T},$$
where $L_{s}^{c}$ and $L_{s}^{d}$ denote the continuous part and jump part of $L_{s}$ respectively, hence $(B_{T}^{t,\theta})_{T\geq t}$ is a L\'{e}vy stochastic integral.
In the following we will prove that, for each $n=1,2,\ldots,$  any $\varphi\in C_{b,Lip}(\mathbb{R}^{^{d\times n}})$ and $0\leq t_{1},\ldots, t_{n}<\infty,$ the $\mathbb{E}^{G}$ defined in \cite{HP 2009a} can be equivalently defined by
$$\mathbb{E}^{G}[\varphi(X_{t_{1}}, X_{t_{2}}-X_{t_{1}}, \ldots X_{t_{n}}-X_{t_{n-1}})]=\sup_{\theta \in \mathcal {A}_{t,T}^{\Theta}}E_{P}[\varphi(B_{t_{1}}^{0,\theta}, B_{t_{2}}^{t_{1},\theta},\ldots, B_{t_{n}}^{t_{n-1},\theta} )].$$

If for each given $\varphi\in C_{b, Lip}(\mathbb{R}^{d})$, we set  $v(t,x):= \sup _{\theta\in \mathcal {A}_{t, T}^{\Theta}}E_{P}[\varphi(x+ B_{T}^{t, \theta})]$, then we can get the following generalized dynamical programming principle:
\begin{e-proposition}
We have
\begin{equation}\label{dpp}
v(t,x)=\sup _{\theta\in \mathcal {A}_{t, t+h}^{\Theta}}E_{P}[v(t+h, x+ B_{t+h}^{t, \theta})].
\end{equation}
\end{e-proposition}
\proof This is a result treated analogous to \cite{DHP}, therefore we omit the proof here.  $\Box$
\begin{lemma}
 $v$ is bounded by $\sup|\varphi|$. It is  Lipschitz function in $x$ and h$\ddot{o}$lder continuous  in $t$.
\end{lemma}
\proof
Obviously we know that
  $$\sup _{\theta\in \mathcal {A}_{t, t+h}^{\Theta}}E_{P}[v(t+h, x+ B_{t+h}^{t, \theta})-v(t+h,x)]=v(t,x)-v(t+h,x),$$
since $v$  is Lipschitz function in $x$, the absolute value of left hand is bounded by
  $$C \sup_{\theta\in \mathcal {A}_{t, t+h}^{\Theta}}E_{P}[|B_{t+h}^{t, \theta}|]\leq C_{1}(h+h^{\frac{1}{2}}).$$
We get the result. $\Box$

\begin{theorem}\label{G-intp}
$v$ is a viscosity solution of the integro-partial differential equation:

\begin{equation}\label{visco solu proof}
\begin{array}{rcl}
\dfrac{\partial v}{\partial t}+G(v,Dv, D^{2}v)
 & = &0, \ \  on \ (t,x)\in [0,T)\times \mathbb{R}^{d}, \\[10pt]
 v(T,x)& = & \varphi(x),
\end{array}
\end{equation}
where the function $G$ is given in (\ref{ipde}).
\end{theorem}
\proof
Since $X_{t}$ has a decomposition with $X_{t}=X_{t}^{c}+X_{t}^{d}$ and for each $t\geq0,$ $L_{t}$ has the following L\'{e}vy-It\^{o} decomposition (see e.g. \cite{Apple}, \cite{Bertoin} or \cite{Sato}):
$$L_{t}=L_{t}^{c}+L_{t}^{d}=bt+W_{t}+\int_{|x|<1}x\tilde{N}(t,dx)+\int_{|x|\geq1}x N(t,dx),$$
where $W$ is a Brownian motion, $N$ is an independent Poisson random measure on $\mathbb{R}^{+}\times (\mathbb{R}^{d}\backslash \{0\})$ with levy measure $\nu$, and $\tilde{N}$ is the compensated Poisson random measure. Then we can write $B_{T}^{t,\theta}$ in
the following form:
 $$B_{T}^{t,\theta}=\int_{t}^{T}{}^{1}\!\theta^{c} b ds+\int_{t}^{T}{}^{2}\!\theta^{c}dW_{s}+\int_{t}^{T}\int_{|x|<1}\theta_{s}^{d} x\tilde{N}(ds,dx)+\int_{t}^{T}\int_{|x|\geq1}\theta_{s}^{d} xN(ds,dx).$$

Firstly, we will suppose that $\nu(\mathbb{R}^{d})<\infty$ which means that $L_{t}$ has finite activity (i.e. it has a finite number of jumps in any finite period of time).\\
Let $B_{T}^{t,\theta,c}$ and $B_{T}^{t,\theta,d}$ be the continuous part and the discontinuous part of $B_{T}^{t,\theta}$ defined respectively by
$$
B_{T}^{t,\theta,c}
  = \int_{t}^{T}{}^{1}\!\theta_{s}^{c}b ds+\int_{t}^{T}{}^{2}\!\theta_{s}^{c}dW_{s}$$
  and $$
B_{T}^{t,\theta,d} = \int_{t}^{T}\int_{|x|<1}\theta_{s}^{d} x\tilde{N}(ds,dx)+\int_{t}^{T}\int_{|x|\geq1}\theta_{s}^{d} xN(ds,dx). $$

Let $\psi\in C_{b}^{2,3}((0,T)\times\mathbb{R}^{d})$ be such that $\psi\geq v$ and, for a fixed $(t,x)\in (0,T)\times\mathbb{R}^{d}, \psi(t,x)=v(t,x).$
By It\^{o}'s formula for L\'{e}vy-type stochastic integral, it follows that
$$\begin{array}{rl}
& \psi(t+h, x+B_{t+h}^{t, \theta})-\psi(t,x)\\[10pt]
 = & \int_{t}^{t+h}\dfrac{\partial \psi}{\partial s}(s, x+B_{s}^{t, \theta})ds  +\dfrac{1}{2}\int_{t}^{t+h}tr[{}^{2}\!\theta_{s}^{c}{{}^{2}\!\theta_{s}^{c}}^{T}D^{2}\psi](s, x+B_{s-}^{t, \theta})ds+\int_{t}^{t+h}\langle D\psi(s, x+B_{s-}^{t, \theta}), dB_{s}^{t, \theta,c}\rangle\\[10pt]
 &+\int_{t}^{t+h} \int_{|z|\geq1}\big[\psi(s,x+B_{s-}^{t, \theta}+\theta_{s}^{d}z)-\psi(s, x+ B_{s-}^{t, \theta})\big]N(ds,dz)\\[10pt]
 & +\int_{t}^{t+h} \int_{|z|<1}\big[\psi(s,x+B_{s-}^{t, \theta}+\theta_{s}^{d}z)-\psi(s, x+ B_{s-}^{t, \theta})\big]\tilde{N}(ds,dz)\\[10pt]
 & +\int_{t}^{t+h} \int_{|z|<1}\big[\psi(s,x+B_{s-}^{t, \theta}+\theta_{s}^{d}z)-\psi(s, x+ B_{s-}^{t, \theta})-\langle\theta_{s}^{d}z, D\psi(s, x+B_{s-}^{t, \theta})\rangle\big]ds \nu(dz)\\[10pt]
 \triangleq & \int_{t}^{t+h}I_{1}ds+\int_{t}^{t+h}\langle D\psi(s, x+B_{s-}^{t, \theta}), {}^{2}\!\theta_{s}^{c}dW_{s}\rangle\\[10pt]
 &+\int_{t}^{t+h} \int_{|z|\geq1}I_{2}N(ds,dz)+\int_{t}^{t+h} \int_{|z|<1}I_{2}\tilde{N}(ds,dz)-\int_{t}^{t+h} \int_{|z|<1}I_{3}ds \nu(dz).\end{array}$$
Obviously, we have  $E_{P}\big[\int_{t}^{t+h}\langle D\psi(s, x+B_{s-}^{t, \theta}), {}^{2}\!\theta_{s}^{c}dW_{s}\rangle\big]=0.$\\
 The uniformly Lipschitz continuity of $\big(\dfrac{\partial\psi}{\partial s}+\langle D\psi,{}^{1}\!\theta_{s}^{c}b\rangle+\dfrac{1}{2}tr[{{}^{2}\!\theta_{s}^{c} {}^{2}\!\theta_{s}^{c}}^{T}D^{2}\psi]\big)(s,y)$ in $(s,y)$
  yields
$$\begin{array}{rl}
E_{P}[I_{1}] & =E_{P}\big[\dfrac{\partial \psi}{\partial s}+\langle D\psi,{}^{1}\!\theta_{s}^{c}b\rangle+\dfrac{1}{2}tr[{}^{2}\!\theta_{s}^{c} {}^{2}\!{\theta_{s}^{c}}^{T}D^{2}\psi](s, x+B_{s-}^{t, \theta})\big]\\[10pt]
 & \leq E_{P}\big[\dfrac{\partial \psi}{\partial s}+\langle D\psi, {}^{1}\!\theta_{s}^{c}b\rangle+\dfrac{1}{2}tr[{{}^{2}\!\theta_{s}^{c} {}^{2}\!\theta_{s}^{c}}^{T}D^{2}\psi]\big](t, x)+C_{1}(h+h^{1/2}), \end{array}$$

\begin{equation}\label{ie1}
\begin{array}{rl}
 & E_{P}\big[\int_{t}^{t+h} \int_{|z|\geq1}I_{2}N(ds,dz)\big]\\[10pt]
 \leq &E_{P}\big[\int_{t}^{t+h} \int_{|z|\geq1}\big(\psi(t,x+\theta_{s}^{d}z)-\psi(t,x)+C_{1}(h+h^{1/2})\big)N(ds,dz)\big] \end{array}
 \end{equation}
and
\begin{equation}\label{ie2}
\begin{array}{rl}
 & E_{P}\big[\int_{t}^{t+h} \int_{|z|<1}\langle\theta_{s}^{d}z, D\psi(s, x+B_{s-}^{t, \theta}\rangle ds \nu(dz)\big]\\[10pt]
 \leq & E_{P}\big[\int_{t}^{t+h} \int_{|z|<1}(\langle\theta_{s}^{d}z, D\psi(t, x)\rangle+C_{1}(h+h^{1/2})|\theta_{s}^{d}z|)ds \nu(dz)\big].\end{array}
 \end{equation}
Since $\theta_{s}$ is a predictable process, then it is independent of $N(ds,dz)$ and $\tilde{N}(ds,dz)$, henceforth we have
$$E_{P}\big[\int_{t}^{t+h} \int_{|z|<1}I_{2}\tilde{N}(ds,dz)\big]=0$$
and
$$\begin{array}{rl}
  & E_{P}\big[\int_{t}^{t+h} \int_{|z|\geq1}\big(\psi(t,x+\theta_{s}^{d}z)-\psi(t,x)\big)N(ds,dz)\big]\\[10pt]
 =& \int_{t}^{t+h} \int_{|z|\geq1} E_{P}\big(\psi(t,x+\theta_{s}^{d}z)-\psi(t,x)\big)\nu(dz)ds. \end{array}$$
Thus, by the assumption on $\nu$ and the inequalities (\ref{ie1}) and (\ref{ie2}), we have
$$\begin{array}{rl}
 & E_{P}\big[\int_{t}^{t+h} \int_{|z|\geq1}I_{2}N(ds,dz)-\int_{t}^{t+h} \int_{|z|<1}I_{3}ds \nu(dz)\big]\\[10pt]
 \leq& \int_{t}^{t+h} \int_{\mathbb{R}^{d}\backslash\{0\}} E_{P}\big(\psi(t,x+\theta_{s}^{d}z)-\psi(t,x)\big)\nu(dz)ds-E_{P}\big[\int_{t}^{t+h} \int_{|z|<1}\langle D\psi(t, x),\theta_{s}^{d}z \rangle \nu(dz) ds\big]\\[10pt]
 & +C_{2}(h^{2}+h^{3/2})\\[10pt]
 \leq& \int_{t}^{t+h} \int_{\mathbb{R}^{d}\backslash\{0\}} \int_{\mathbb{R}^{d}\backslash\{0\}}\big(\psi(t,x+z')-\psi(t,x)\big)F_{\theta_{s}^{d}}(d(\dfrac{z'}{z}))\nu(dz)ds
 +E_{P}\big[\int_{t}^{t+h} \langle D\psi(t, x),\theta_{s}^{d}b' \rangle ds\big]\\[10pt]
 &+C_{2}(h^{2}+h^{3/2}))\\[10pt]
 =& h\int_{\mathbb{R}^{d}\backslash\{0\}}\big(\psi(t,x+z')-\psi(t,x)\big)\nu'(dz')
 +E_{P}\big[\int_{t}^{t+h} \langle D\psi(t, x),\theta_{s}^{d}b' \rangle ds\big]\\[10pt]
 &+C_{2}(h^{2}+h^{3/2}),\\[10pt]
 \end{array}$$
where $F_{\theta_{s}^{d}}$ is the probability distribution function
of $\theta_{s}^{d}$, and we denote
$\nu'(dz')=\int_{\mathbb{R}^{d}\backslash\{0\}}F_{\theta_{s}^{d}}(d(\dfrac{z'}{z}))\nu(dz)$,
and $b'=-\int_{|z|<1}z\nu(dz)$.\\
Hence, from the dynamic programming principle (\ref{dpp}) it follows that
$$\begin{array}{rcl}
0& = & \sup_{\theta\in \mathcal {A}_{t,t+h}^{\Theta}}E_{P}[v(t+h, x+B_{t+h}^{t, \theta})-v(t,x)] \\[10pt]
& \leq & \sup_{\theta\in \mathcal {A}_{t,t+h}^{\Theta}}E_{P}[\psi(t+h, x+B_{t+h}^{t, \theta})-\psi(t,x)] \\[10pt]
& \leq & \sup_{\theta\in \mathcal {A}_{t,t+h}^{\Theta}}\int_{t}^{t+h}E_{P}\big[\dfrac{\partial \psi}{\partial s}+\langle D\psi,{}^{1}\!\theta_{s}^{c}b+\theta_{s}^{d}b')\rangle+\dfrac{1}{2}tr[{}^{2}\!\theta_{s}^{c} {}^{2}\!{\theta_{s}^{c}}^{T}D^{2}\psi]\big](t, x)ds+C_{1}(h^{2}+h^{3/2})      \\[10pt]
& & + \sup_{\nu'\in \mathcal {V}}h\int_{\mathbb{R}^{d}\backslash\{0\}}\big(\psi(t,x+z')-\psi(t,x)\big)\nu'(dz')+C_{2}(h^{2}+h^{3/2})    \\[10pt]
& \leq & h\sup _{(\nu',q,\gamma)\in \mathcal{U}}\{\int_{\mathbb{R}^{d}\backslash \{0\}}(\psi(t,x+z)-\psi(t,x))\nu'(dz)+\langle D\psi(t,x),q\rangle+ \frac{1}{2}tr[\gamma\gamma^{T}D^{2}\psi(t,x)]\} \\[10pt]
& & +h \dfrac{\partial\psi}{\partial s}(t,x)+(C_{1}+C_{2})(h^{2}+h^{3/2}).  \\[10pt]\end{array}$$
Consequently,
$$
\frac{\partial\psi}{\partial s} (t,x)+\sup _{(\nu',q,\gamma)\in \mathcal{U}}\{\int_{\mathbb{R}^{d}\backslash \{0\}} (\psi(t,x+z)-\psi(t,x))\nu'(dz)+\langle D\psi(t,x),q\rangle+ \frac{1}{2}tr[\gamma\gamma^{T}D^{2}\psi(t,x)]\}\geq0.$$
By the definition, $v$ is a viscosity subsolution.
Now we will prove that this conclusion is also true if we remove the condition $\nu(\mathbb{R}^{d})<\infty$. \\
Denote
$${}^{\varepsilon}\!{L_{t}}=bt+W_{t}+\int_{\varepsilon\leq|x|<1}x\tilde{N}(t,dx)+\int_{|x|\geq1}x N(t,dx)$$
and $$^{\varepsilon}\!B_{T}^{t}=B_{T}^{t,\theta,c}+{}^{\varepsilon}\!B_{T}^{t,\theta,d},$$
where
$${}^{\varepsilon}\!B_{T}^{t,\theta,d} = \int_{t}^{T}\int_{\varepsilon\leq|x|<1}\theta_{s}^{d} x\tilde{N}(ds,dx)+\int_{t}^{T}\int_{|x|\geq1}\theta_{s}^{d} xN(ds,dx),$$
note that each $({}^{\varepsilon}\!{L_{t}})_{t\geq0}$ is a a compound poisson process, hence is a finite-activity process.
$$\begin{array}{rcl}
0& = & \sup_{\theta\in \mathcal {A}_{t,t+h}^{\Theta}}E_{P}[v(t+h, x+B_{t+h}^{t, \theta})-v(t,x)] \\[10pt]
 & \leq& \sup_{\theta\in \mathcal {A}_{t,t+h}^{\Theta}}E_{P}[v(t+h, x+{}^{\varepsilon}\!B_{t+h}^{t, \theta})-v(t,x)]+ C\sup_{\theta\in \mathcal {A}_{t,t+h}^{\Theta}}E_{P}[|\int_{t}^{t+h}\int_{|x|<\varepsilon}\theta_{s}^{d} z\tilde{N}(ds,dz)|]\\[10pt]
& \leq & \sup_{\theta\in \mathcal {A}_{t,t+h}^{\Theta}}E_{P}[\psi(t+h, x+{}^{\varepsilon}\!B_{t+h}^{t, \theta})-\psi(t,x)]+C\sup _{\theta\in \mathcal {A}_{t,t+h}^{\Theta}}\int^{t+h}_{t}\int_{|x|<\varepsilon}E_{P}[|\theta_{s}^{d}z|]\nu(dz)ds \\[10pt]
& \leq & \sup_{\theta\in \mathcal {A}_{t,t+h}^{\Theta}}\int_{t}^{t+h}E_{P}\big[\dfrac{\partial \psi}{\partial s}+\langle D\psi,{}^{1}\!\theta_{s}^{c}b+\theta_{s}^{d}{}^{\varepsilon}\!b')\rangle+\dfrac{1}{2}tr[{}^{2}\!\theta_{s}^{c} {}^{2}\!{\theta_{s}^{c}}^{T}D^{2}\psi]\big](t, x)ds+C_{1}(h^{2}+h^{3/2})      \\[10pt]
& & +\sup_{\theta\in \mathcal {A}_{t,t+h}^{\Theta}} \int_{t}^{t+h} \int_{|z|\geq\varepsilon} E_{P}\big(\psi(t,x+\theta_{s}^{d}z)-\psi(t,x)\big)\nu(dz)ds+C_{2}(h^{2}+h^{3/2})+C_{3}h\int_{|z|<\varepsilon}|z|\nu(dz),\\[10pt]
\end{array}$$
where ${}^{\varepsilon}\!b'=\int_{\varepsilon\leq|z|}|z|\nu(dz)$. Here note that $C_{2}$ is a finite number for each fixed $\varepsilon$. Since $h>0$ is arbitrary small, and $\int_{\mathbb{R}^{d}\backslash
\{0\}}|z|\nu(dz)<\infty$ implies that
 $\lim_{\varepsilon\downarrow0}\int_{|z|<\varepsilon}|z|\nu(dz)=0$.
Then we divide both sides of last inequality by $h$,
 and let firstly $h$ and then $\varepsilon$ go to zero, we could finally obtain the result.

 Similarly, we can prove that $v$ is also a supersolution.
The proof is complete now. $\Box$

We know that $u(t,x):=v(T-t,x)$, then $u$ is a viscosity solution of $\dfrac{\partial u}{\partial t}-G(u,Du, D^{2}u) = 0$, with initial condition $u(0,x)=\varphi(x)$.
From the uniqueness of the viscosity solution of integro-pde and Theorem \ref{G-intp}, we get the following
proposition:
\begin{e-proposition}
$$\begin{array}{rl}
\mathbb{E}^{G}[\varphi(X_{t_{1}}, X_{t_{2}}-X_{t_{1}}, \ldots X_{t_{n}}-X_{t_{n-1}})]&=\sup_{\theta \in \mathcal {A}_{0,\infty}^{\Theta}}E_{P}[\varphi(B_{t_{1}}^{0,\theta}, B_{t_{2}}^{t_{1},\theta},\ldots, B_{t_{n}}^{t_{n-1},\theta} )]\\[10pt]
 &=\sup_{\theta \in \mathcal {A}_{0,\infty}^{\Theta}}E_{P_{\theta}}[\varphi(B_{t_{1}}^{0}, B_{t_{2}}^{t_{1}},\ldots, B_{t_{n}}^{t_{n-1}} )], \end{array}$$
where $P_{\theta}$ is the law of the process $B_{t}^{0,\theta}$, $t\geq0$, for $\theta \in \mathcal {A}_{0,\infty}^{\Theta}$.
\end{e-proposition}

\begin{e-proposition} The family of probability measures $\{P_{\theta}, \theta \in \mathcal {A}^{\Theta}_{0,\infty}\}$ on $\Omega=D_0([0,\infty),\mathbb{R}^{d})$ is tight.
\end{e-proposition}
\proof Since for any $ s,u,t \ \  0\leq s\leq u\leq t <\infty,$ there exists $C'>0$, such that  $\mathbb{E}^{G}[|X^{c}_{t}-X^{c}_{s}|^{4}]\leq C'|t-s|^{2},$ (see. Proposition 49 of \cite{DHP})
and

$$\begin{array}{rl}
\mathbb{E}^{G}[|X^{d}_{t}-X^{d}_{u}|\cdot|X^{d}_{u}-X^{d}_{s}|]&=\sup_{\theta \in \mathcal {A}_{t,T}^{\Theta}}E_{P_{\theta}}[|B_{t}^{u,\theta,d}|\cdot|B_{u}^{s,\theta,d}|]\\[10pt]
 &\leq C'|t-s|^{2}. \end{array}$$
Therefore the statement follows from Corollary \ref{kca} and  Theorem \ref{klc}. $\Box$
\begin{example}
Let $(L_{t})_{t\geq0}$ be a homogeneous Poisson process with intensity $\lambda=1$ denoted by $(N_{t})_{t\geq0}$, then we take especially a collection of $\mathbb{F}$-predictable process $\mathcal {A}_{t,T}^{\Theta^{d}}$, where $\Theta^{d}=\{0,1\}$ and for each
$t>0$, $\theta_{t}^{d}$ follows a Bernoulli distribution with success probability $p$, where $p\in [\lambda,1]$.
If $(X_{t})_{t\geq0}$ is the 1-dimensional $G$-Poisson process \cite{HP 2009b} defined  by the following
equation:
$$\partial_{t}u(t,x)-G_{\lambda}(u(t,x+1)-u(t,x))=0, \ \  u(0,x)=\varphi(x),$$
where $G_{\lambda}(a)=a^{+}-\lambda a^{-}, \lambda\in[0,1].$ It is easy to check that, for each $n\in \mathbb{N}$,  any $\varphi\in C_{b,Lip}(\mathbb{R}^{^{d\times n}})$ and $0\leq t_{1},\ldots, t_{n}<\infty,$
$$\mathbb{E}^{G}[\varphi(X_{t_{1}}, X_{t_{2}}- X_{t_{1}}, \ldots X_{t_{n}}-X_{t_{n-1}})]=\sup_{\theta^{d} \in \mathcal {A}_{0,\infty}^{\Theta^{d}}}E_{P}[\varphi(\int_{t_{1}}^{0}\theta_{s}^{d}dN_{s}, \int_{t_{2}}^{t_{1}}\theta_{s}^{d}dN_{s},\ldots,\int_{t_{n}}^{t_{n-1}}\theta_{s}^{d}dN_{s})].$$
\end{example}



\section{Characterization of $L_{ip}(\Omega)$}
In \cite{DHP}, we know that if \
$\Omega_{T}=C_0([0,T],\mathbb{R}^{d})$  for $T>0,$
(resp. $\Omega=C_0([0,\infty),\mathbb{R}^{d})$), then \ for  $t\leq T$,
 $$L_{ip}(\Omega_{t})\subseteq L_{ip}(\Omega_{T})\subset C_{b}(\Omega_{T}) \ \ (resp. \ L_{ip}(\Omega)\subset C_{b}(\Omega)),$$
and that an element $Y$ of the space $L_{G}^{p}(\Omega)$ is a quasi-continuous function $Y=Y(\omega)$ defined on $\Omega$. $L_{G}^{p}(\Omega)$ is also proved to be identified with  the space $\mathscr{L}^{p}$ that introduced in \cite{DM}.

 But if $\Omega=D_0([0,T]),\mathbb{R}^{d})$, this relationship of inclusion between $ L_{ip}(\Omega_{T})$ and $C_{b}(\Omega_{T})$ is no longer true. In fact, for any fixed $t\in[0,T]$, the variable $B_{t}(\omega)=\omega_{t}$ (or $\omega(t)$) on $(\Omega, \mathscr{B}(\Omega))$
 $$\begin{array}{cccc}
B_{t}: & \Omega & \mapsto & \mathbb{R}\\
 & \omega & \mapsto & \omega_{t}\end{array}$$
is continuous in $\omega$ if and only if $\omega$ is continuous at $t$ (see Section 12 of \cite{Bill}), thus $L_{ip}(\Omega_{T})$ does not belong to $C_{b}(\Omega_{T})$.


\begin{proposition}
Let $\Omega_{T}=D_0([0,T],\mathbb{R})$ be equipped with the the Skorohod metric $d_{s}$, then $B_{t}$ is either upper semi-continuous (in short, u.s.c.) or lower semi-continuous (in short, l.s.c.) at each point $\omega \in \Omega_{T}.$
\end{proposition}
\proof For any $\omega_{\circ}\in \Omega_{T},$
$\omega_{\circ}$ is either u.s.c. or l.s.c. at each fixed $t\in[0,T]$. Without loss of generality ,  we will only consider the
u.s.c. case. Then for any $\varepsilon>0$, there exists $\delta >0,$ such
that for any $s\in[0,T]$ satisfying that $|s-t|<\delta$, we have $\omega_{\circ}(s)< \omega_{\circ}(t)+\varepsilon/2$. Now we can choose $\eta=\min(\varepsilon/2, \delta)>0$, for any
$\omega$ such that $d_{s}(\omega,\omega_{\circ})<\eta$, there exists some $\lambda$ which is a strictly increasing, continuous mapping from [0,T] onto itself with $\lambda 0=0$ and $\lambda T=T$, then we have $\sup _{t}|\omega(t)-\omega_{\circ}(\lambda(t))|<\eta$, and $\sup _{t}|\lambda(t)-t|<\eta$. Hence we get
$$\omega(t)-\omega_{\circ}(t)=\omega(t)-\omega_{\circ}(\lambda(t))+\omega_{\circ}(\lambda(t))-\omega_{\circ}(t)
<\eta+\varepsilon/2<\varepsilon.$$
Thus the proof is complete. $\Box$\\
\begin{proposition}\label{quasic}
For each $X\in L_{ip}(\Omega)$ and $\varepsilon>0$, there exists $Y\in C_{b}(\Omega)$ such that
$\mathbb{E}^{G}[|X-Y|]<\varepsilon$.
\end{proposition}
\proof
Without loss of generality, we suppose that $\Omega=D_0([0,\infty), \mathbb{R}).$ Let each $X\in L_{ip}(\Omega)$ be with the form $X=\varphi(B_{t_{1}},B_{t_{2}},\cdots B_{t_{m}})$, where $\varphi\in C_{b, Lip}(\mathbb{R}^{m})$ and $0\leq t_{1}\leq t_{2}\leq\cdots\leq t_{m}< \infty$. For each $i=1,2,\cdots, m,$ let $h^{i}_{n}(\omega)=n\int_{t_{i}}^{t_{i}+\frac{1}{n}}B_{s}(\omega)ds,$ then $h^{i}_{n}$ is continuous in the Skorohod topology(see\cite{Bill}). In fact, if $\omega_{k}\rightarrow \omega$  in the Skorohod topology, then $\omega_{k}(s)\rightarrow \omega(s)$ for continuity points $s$ and hence for points $s$ outside a set of Lebesgue measure 0; since $\omega_{k}$ are uniformly bounded, we have $\lim_{k} h^{i}_{n}(\omega_{k})\rightarrow h^{i}_{n}(\omega).$ By the right continuity of $\omega$, $h^{i}_{n}(\omega)\rightarrow B_{t_{i}}(\omega)$ as $n\rightarrow\infty$. Then it follows that

$$\begin{array}{ccl}
 \mathbb{E}^{G}[|\varphi(B_{t_{1}},B_{t_{2}},\cdots B_{t_{m}})-\varphi(h_{n}^{1},h_{n}^{2},\cdots h_{n}^{m})|]
 & \leq & \sum_{i=1}^{m}\mu\mathbb{E}^{G}[|B_{t_{i}}-h_{n}^{i}|]\\[10pt]
 & = & \sum_{i=1}^{m}\mu\mathbb{E}^{G}[|B_{t_{i}}-n\int_{t_{i}}^{t_{i}+\frac{1}{n}}B_{s}ds|]\\[10pt]
 & = & \sum_{i=1}^{m}\mu\mathbb{E}^{G}[|n\int_{t_{i}}^{t_{i}+\frac{1}{n}}(B_{t_{i}}-B_{s})ds|]\\[10pt]
 & \leq& \sum_{i=1}^{m}\mu n\int_{t_{i}}^{t_{i}+\frac{1}{n}}\mathbb{E}^{G}[|B_{t_{i}}-B_{s}|]ds\\[10pt]
 & \leq & \sum_{i=1}^{m}\mu n\int_{t_{i}}^{t_{i}+\frac{1}{n}}C(|t_{i}-s|^{1/2}+|t_{i}-s|)ds\end{array},$$
where $\mu$ is the Lipschitz constant of $\varphi$. Hence for each positive $\varepsilon$, we can choose some $n_{0}>0$ and set $Y=\varphi(h_{n_{0}}^{1},h_{n_{0}}^{2},\cdots h_{n_{0}}^{m})$ such that $\mathbb{E}^{G}[|X-Y|]<\varepsilon$. $\Box$
\begin{remark}
This proposition implies that $L_{ip}(\Omega)\subseteq \mathbb{L}_{c}^{1}$, hence $L_{G}^{1}(\Omega)\subseteq \mathbb{L}_{c}^{1}.$
\end{remark}

As shown in previous sections, we use two different methods to prove that $\mathbb{E}^{G}$ is an upper expectation associated to a weakly relatively compact family $\mathscr{P}_{1}$ and now denote by $\mathscr{P}=\overline{\mathscr{P}_{1}}$ the closure of $\mathscr{P}_{1}$ under the topology  of  weak convergence, then $\mathscr{P}$ is weakly compact. For each $X\in L^{0}(\Omega)$ such that $E_{P}[X]$ exists for each $P\in\mathscr{P}$, we set
$$\mathbb{E}^{\mathscr{P}}[X]=\sup_{P\in \mathscr{P}}E_{P}[X]$$
and
$$\mathbb{E}^{\mathscr{P}_{1}}[X]=\sup_{P\in \mathscr{P}_{1}}E_{P}[X].$$
Proposition \ref{quasic} together with Proposition \ref{wc} yields the following theorem.
\begin{theorem}
For each $X\in L_{G}^{1}(\Omega)$, we have $\mathbb{E}^{G}[X]=\mathbb{E}^{\mathscr{P}}[X]=\mathbb{E}^{\mathscr{P}_{1}}[X]$.
\end{theorem}
{\bf Acknowledgements}
\\
\\
The author is grateful to professor Shige Peng for his encouragement
and helpful discussions.

\selectlanguage{english}
\appendix
\section{}
\begin{lemma}\label{kl} (See \cite{Bill})
Let $\{P_{n}\}$ be a sequence of probability measures on a measurable space $({D_0([0,T],\mathbb{R})},\mathcal {D})$ (resp. $({D_0([0,\infty),\mathbb{R})},\mathcal {D}_{\infty})$), then $\{P_{n}\}$ is tight if and only if these two conditions hold:\\
(i) For each t in a set $\mathcal {T}$ dense in [0,T] and contains T (resp. dense in $[0,\infty)$),
$$\lim_{a\rightarrow\infty}\lim\sup_{n}P_{n}[\omega:|\omega_{t}|\geq a]=0.$$
(ii) For each positive $\varepsilon$ (resp. for each positive $\varepsilon$ and T),
$$
\begin{cases}\begin{array}{l}
\lim_{\delta\rightarrow0}\lim\sup_{n}P_{n}[\omega:V''_{T}(\omega,\delta)|\geq \varepsilon]=0.\\[10pt]
\lim_{\delta\rightarrow0}\lim\sup_{n}P_{n}[\omega:|\omega_{\delta}|\geq \varepsilon]=0.\\[10pt]
\lim_{\delta\rightarrow0}\lim\sup_{n}P_{n}[\omega:|\omega_{T-}-\omega_{T-\delta}|\geq \varepsilon]=0.\end{array}\end{cases}$$
where $V''_{T}(\omega,\delta)=\sup _{s\leq u\leq t, \ t-s\leq\delta}
\{|\omega_{u}-\omega_{s}|\wedge|\omega_{t}-\omega_{u}|\},$ the supremum extending over
all triples $s, u, t$ in [0,T] satisfying the constraints.
\end{lemma}
\begin{corollary}\label{kca}
Let X and Y be two  c\`{a}dl\`{a}g stochastic processes. Assume that the distributions of X and Y respectively under the probability measures $(P_{n}, n\in \mathbb{N})$ denoted respectively by $\big(P_{n}(X\in\cdot), n\in \mathbb{N}\big)$ and $\big(P_{n}(Y\in\cdot), n\in \mathbb{N}\big)$ are both tight,
then $\big(P_{n}(X+Y\in\cdot), n\in \mathbb{N}\big)$ is tight.
\end{corollary}
\proof
Apply the triangle inequality in order to check the conditions in the preceding lemma. $\Box$
\begin{theorem}\label{klc}(Kolmogorov-Chentsov's criterion for weak relative compactness)
Let $\mathscr{P}$ be any subset of the collection of all probability
measures on ${D_0([0,T],\mathbb{R})}$ and $\mathbb{E}$ the upper
expectation related to $\mathscr{P}$.
If the following conditions are satisfied:\\
(i) $\exists a >0$, such that $\mathbb{E}[|\omega_{t}-\omega_{s}|]\leq C|t-s|^{a}$, $\forall t,s \in[0,T] $;\\
(ii)
$\mathbb{E}[|\omega_{t}-\omega_{u}|^{p}|\omega_{u}-\omega_{s}|^{q}]\leq
C|t-s|^{1+r},$ for some $C, r >0, p, q\geq0$ with $p+q >0,$  and all
$0\leq s\leq u\leq t \leq T$, then $\mathscr{P}$ is relatively
compact.
\end{theorem}
\proof
Let $\{P_{n}\}$ be any sequence in $\mathscr{P}$. We check the conditions of the previous lemma. Obviously,
condition (i) implies $\mathbb{E}[|\omega_{\delta}|]\leq C\delta^{a}$ and \  $\mathbb{E}[|\omega_{T-}-\omega_{T-\delta}|]\leq C\delta^{a}$.
 By condition(ii), $\forall\eta>0, \alpha\in(0, \dfrac{r}{p+q}), \ \exists K>0,$ for every $n$, we have
$$P_{n}\big[|\omega_{t}-\omega_{u}|\geq K|t-s|^{\alpha}, |\omega_{u}-\omega_{s}|\geq K|t-s|^{\alpha},\  \forall 0\leq s\leq u\leq t\leq T]\big]\leq\eta,$$
This obviously implies $\lim\sup_{n} P_{n}[V_{T}''(\omega,\delta)\geq \varepsilon]\leq\eta$ with  $\delta=\big|\dfrac{\varepsilon}{K}\big|^{1/\alpha}$.
Hence by Lemma \ref{kl}, we get the conclusion. $\Box$\\

\end{document}